\documentclass{article}
\usepackage[left=3cm,right=3cm,top=3.5cm,bottom=3.5cm]{geometry}
\PassOptionsToPackage{dvipsnames,table}{xcolor}
\usepackage[toc,page]{appendix}
\usepackage[utf8]{inputenc}
\usepackage[T1]{fontenc}
\usepackage{tikz-cd}
\usepackage{tikz}
\usepackage{nextpage}
\usepackage{circledsteps}
\usepackage{amsfonts}
\usepackage{amsmath}
\usepackage{amssymb}
\usepackage{amsthm}
\usepackage{mathtools}
\usepackage[english]{babel}
\usepackage{enumitem}
\usepackage{comment}
\usepackage{makecell}
\usepackage{mathrsfs}
\usepackage{mathdots}\usepackage[Sonny]{fncychap}
\usepackage{adjustbox}
\usepackage{thmtools}
\usepackage{textcomp}
\usepackage{float}
\usepackage[hidelinks]{hyperref}
\usepackage{cleveref}
\usepackage{arydshln}

\DeclareMathOperator{\Div}{Div}

\DeclareMathOperator{\Spec}{Spec}

\DeclareMathOperator{\RG}{R\Gamma}
\DeclareMathOperator{\HH}{H}

\DeclareMathOperator{\id}{id}

\DeclareMathOperator{\TC}{TC}
\DeclareMathOperator{\Supp}{Supp}
\DeclareMathOperator{\Aut}{Aut}

\DeclareMathOperator{\im}{im}
\DeclareMathOperator{\Hom}{Hom}

\DeclareMathOperator{\Gal}{Gal}
\DeclareMathOperator{\Pic}{Pic}

\newcommand{\LL}{\mathcal{L}}

\newcommand{\R}{{\rm R}}

\newcommand{\ZZ}{\mathbb{Z}}
\newcommand{\PP}{\mathbb{P}}

\renewcommand{\SS}{\Gamma}

\newcommand{\GG}{\mathfrak{G}}
\newcommand{\F}{\mathcal{F}}

\newcommand{\tX}{{\tilde{X}}}
\newcommand{\tY}{{\tilde{Y}}}

\newcommand{\Xet}{X_{\rm \acute{e}t}}

\newcommand{\Xfet}{X_{\rm f\acute{e}t}}
\newcommand{\Xlet}{X_{\rm \ell\acute{e}t}}
\newcommand{\red}{{\rm red}}
  
\newcommand{\nt}{^{\langle n\rangle}}
\newcommand{\nnt}{^{\langle n\rangle\langle n\rangle}}
\newcommand{\nnnt}{^{\langle n\rangle\langle n\rangle\langle n\rangle}}
\newcommand{\ns}{^{[n]}}
\newcommand{\nns}{^{[n][n]}}
\newcommand{\lt}{^{\langle \ell\rangle}}

\newcommand{\Xnt}{{X\nt}}

\newcommand{\kpf}{k^{\rm pf}}
\newcommand{\Xpf}{X^{\rm pf}}

\newtheorem{Definition}{Definition}[section]
\newtheorem{Remark}[Definition]{Remark}
\newtheorem{theorem}[Definition]{Theorem}
\newtheorem{prop}[Definition]{Proposition}
\newtheorem{Lemma}[Definition]{Lemma}

\newtheorem{Corollary}[Definition]{Corollary}

\usepackage{xpatch}
\makeatletter
\xpatchcmd{\@thm}{\thm@headpunct{.}}{\thm@headpunct{}}{}{}
\makeatother

\title{Curves are algebraic $K(\pi,1)$: theoretical and practical aspects}

\author{Christophe Levrat \footnote{INRIA Saclay, 1 rue Honoré d'Estienne d'Orves, 91120 Palaiseau, France. email: christophe.levrat@math.cnrs.fr}}

\date{}

\begin{document}

\maketitle

\begin{abstract}
We prove that any geometrically connected curve $X$ over a field $k$ is an algebraic $K(\pi,1)$, as soon as its geometric irreducible components have nonzero genus. This means that the cohomology of any locally constant constructible étale sheaf of $\ZZ/n\ZZ$-modules, with $n$ invertible in $k$, is canonically isomorphic to the cohomology of its corresponding $\pi_1(X)$-module. To this end, we explicitly construct some Galois coverings of $X$ corresponding to Galois coverings of the normalisation of its irreducible components. When $k$ is finite or algebraically closed, we precisely describe finite quotients of $\pi_1(X)$ that allow to compute the cohomology groups of the sheaf, and give explicit descriptions of the cup products $\HH^1\times\HH^1\to \HH^2$ and $\HH^1\times \HH^2\to \HH^3$ in étale cohomology in terms of finite group cohomology.
\end{abstract}

In this article, the word \textit{curve} refers to a separated one-dimensional scheme of finite type over a field. In particular, a curve is not necessarily reduced, connected or smooth. The \textit{geometric genus} of an irreducible curve $X$ over an algebraically closed field is that of the associated reduced scheme $X_\red$ as defined in \cite[0BYF]{stacks}. The word \textit{covering} refers to a finite flat surjective map of schemes. Cohomology groups are always either étale or profinite group cohomology.

\section{Introduction}

Cohomology theories come in various flavours: while some, like group cohomology, can be defined in a very straightforward manner by constructing an explicit cochain complex, others, like étale cohomology, have no easier known definition than the one using the abstract machinery of derived functors. The mere computability (in the computer theoretic sense) of étale cohomology groups was only proven a few years ago by Madore and Orgogozo \cite[Th. 0.1]{mo}, using fibrations by affine curves. Their proof relies on the fact that affine curves belong to a specific class of connected schemes, called $K(\pi,1)$ schemes, for which the étale cohomology of a locally constant sheaf with invertible torsion is isomorphic to the Galois cohomology of the associated $\pi_1$-module.

Constructions of $K(\pi,1)$ neighbourhoods have been studied since the advent of étale cohomology, when Artin used them to prove a comparison theorem between the étale cohomology of a variety over the complex numbers and the singular cohomology of its analytification \cite[XI, Th. 4.4]{sga43}. Such neighbourhoods were also constructed by Achinger more recently in the context of log smooth schemes, to prove a comparison theorem between étale cohomology and the cohomology of Faltings' topos used in $p$-adic Hodge theory \cite[Th. 9.5]{achinger_compositio}. 

A few years later, Achinger proved using higher ramification theory that any connected affine scheme of positive characteristic is a $K(\pi,1)$ scheme \cite[Th. 1.1]{achinger_Kpi1}. Concerning projective schemes, there are far fewer positive results. It is a classical fact that all smooth connected curves over separably closed fields except $\PP^1$ are $K(\pi,1)$ schemes. In this article, we prove that this is more generally true for any geometrically connected curve over any field, as soon as each of its geometric irreducible components has nonzero geometric genus. This is done in \Cref{sec:kpi1} using trivialising coverings explicitly described in \Cref{sec:Galoiscov}.

\begin{theorem} Let $X$ be a geometrically connected curve over a field $k$. Denote by $\bar k$ an algebraic closure of $k$. If the geometric genus of every irreducible component of $X\times_k \bar k$ is nonzero then $X$ is a $K(\pi,1)$ scheme.
\end{theorem}

Moreover, the automorphism groups of these coverings allow to compute the cohomology of locally constant étale sheaves of $\ZZ/n\ZZ$-modules, as well as the associated cup products. Using technical results on profinite group cohomology presented in \Cref{sec:galoiscoh}, we state explicit results on the computation of such cup products in \Cref{sec:cup}, when the base field $k$ has cohomological dimension 0 or 1 and $\HH^1(k,\ZZ/n\ZZ)$ is finite. This allows to replace étale cohomology not only with profinite group cohomology, but with finite group cohomology, which is easily computable. While the emphasis in the whole article lies on explicit constructions, the computations we present can easily be turned into an algorithm, though it might not be efficient enough to be used in practice.

\section{Galois coverings of algebraic curves}\label{sec:Galoiscov}

In this section, given a curve $X$ over an algebraically closed field as well as a Galois covering of the normalisation of one of its irreducible components, we describe a Galois covering of $X$ with the same Galois group. Similar though not identical results can be found in \cite[§2.3]{das_galois}; we keep the descriptions self-contained, simple and constructive in the specific cases we are interested in. We describe one covering of $X$ which trivialises the étale cohomology groups $\HH^i(X,\ZZ/n\ZZ)$ for $i=1,2$, and explain its contruction explicitly in case $X$ is integral. Finally, we give a similar, though slightly more involved, construction in the case of a curve over a field of cohomological dimension 1. These coverings will be used in the following sections in the proof that curves are $K(\pi,1)$, as well as in the explicit computation of cohomology groups and cup products.

\subsection{Terminology and basic results}\label{subsec:notations}

Given a covering $f\colon Y\to X$ of connected schemes, the notation $\Aut(Y|X)$ refers to the group of automorphisms $\sigma\in\Aut(Y)$ such that $f\circ\sigma=f$. The covering is said to be \textit{Galois} if it is étale and $\Aut(Y|X)$ acts transitively on a geometric fibre of $f$. Transitivity of this action is equivalent to the order of $\Aut(Y|X)$ being exactly the degree of $f$  \cite[03SF]{stacks}.\\
 
Let $X$ be a reduced curve over an algebraically closed field $k$. \begin{itemize}[label=$\diamond$]
\item The notation $\Div(X)$ refers to the group of Weil divisors on $X$.
\item A Weil divisor on $X$ is called simple if it is a sum of distinct closed points of $X$.
\item  Suppose $X$ is smooth. Given simple divisors $D_1,\dots,D_r$ on $X$ with disjoint supports, we will denote by $X_{D_1,\dots,D_r}$ the quotient of $X$ by the equivalence relation whose equivalence classes are the supports of the $D_i$, i.e. the curve defined by the moduli $D_i$ in the sense of \cite[IV.§1.4]{serrega}. The singular points of $X_{D_1,\dots,D_r}$ are the images of the divisors $D_1,\dots,D_r$ of degree at least 2. Such a singularity is called  multicross: the completion of its local ring is of the form $k[[x_1,\dots,x_c]]/(\{x_ix_j\}_{i\neq j})$. 
\item We say that $X$ has at worst multicross singularities if all its singularities are multicross.
\item Suppose $X$ has at worst multicross singularities. Then $X$ may be defined as $Y_{E_1,\dots,E_s}$ where $\nu\colon Y\to X$ is the normalisation map and the $E_j\in\Div(Y)$ are simple. For more details, see for instance \cite[§2.2]{das_galois}. Given simple divisors $D_1,\dots,D_r$ on $X$ with disjoint supports, we will define $X_{D_1,\dots,D_r}$ as follows. For each $i\in \{1\dots r\}$, let $F_i\in \Div(Y)$ be the simple divisor with support $\nu^{-1}(\Supp D_i)$. Let $j_1,\dots,j_t\in \{1\dots s\}$ be those indices $j$ such that the image of $E_{j}$ in $X$ does not lie in the support of any $D_i$. Then we define $X_{D_1,\dots,D_r}$ as $Y_{E_{j_1},\dots,E_{j_t},F_1,\dots,F_r}$. This is a curve with multicross singularities whose normalisation map $Y\to X_{D_1,\dots,D_r}$ factors through $X$.
\item Suppose $X$ has at worst multicross singularities. Consider simple divisors with disjoint support $D_1,\dots,D_r\in \Div(X)$. Pick one of these divisors $D_i=\sum_{j=1}^s P_{i,j}$. Denoting by $Q_i$ the point of $X_{D_1,\dots,D_r}$ whose preimage in $X$ is $D_i$ and by $\TC$ the tangent cones, the construction of $X_{D_1,\dots,D_r}$ yields the following isomorphism: \[ (\TC_{P_{i,1}}(X)\sqcup \dots\sqcup \TC_{P_{i,s}}(X))_{P_{i,1}+\dots+P_{i,s}}\xrightarrow{\sim} \TC_{Q_i}(X_{D_1,\dots,D_r}) .\]
\end{itemize}
More details on curves with multicross singularities can be found in \cite[0C1P]{stacks}.

\subsection{Coverings of irreducible curves with multicross singularities}

\begin{Lemma}\label{lem:seminormcover} Let $X$ be an integral curve with multicross singularities over an algebraically closed field $k$, and $\tilde X$ be its normalisation. For any Galois covering $\tilde Y$ of $\tilde X$, there is a Galois covering $Y\to X$ whose normalisation is $\tilde Y$.
\begin{proof}
Denote by $P_1,\dots,P_r$ the singular points of $X$. For each $P_i$, denote by $P_{i,1},\dots,P_{i,j_i}$ its preimages in $\tilde X$. Denote by $\sigma_1=\id,\sigma_2,\dots,\sigma_d$ the elements of $\Aut(\tY|\tX)$. For each $(i,j)$, take one preimage $P_{i,j}^{(1)}$ of $P_{i,j}$ in $\tY$, and set $P_{i,j}^{(m)}\coloneqq \sigma_m(P_{i,j}^{(1)})$ for all $m\in\{ 1\dots d\}$. For every $i\in \{1\dots r\}$, $m\in \{ 1\dots d\}$, set $D_i^{(m)}=\sum_{j=1}^{j_i}P_{i,j}^{(m)}$. Denote by $Y=\tY_{D_1^{(1)},\dots,D_r^{(d)}}$ the quotient of $\tY$ by the equivalence relation whose equivalence classes are the supports of the divisors $D_i^{(m)}$. Denote by $Q_i^{(m)}$ the image of $D_{i}^{(m)}$ under $\tY\to Y$. The following diagram sums up the notations.
\[
\begin{tikzcd}
Y&& Q_i^{(1)} \arrow[d,no head,dashed] &~ &\arrow[ll,|->] P_{i,1}^{(1)}\arrow[d,no head,dashed] \arrow[r,no head,dashed] & P_{i,j_i}^{(1)}\arrow[d,no head,dashed]&&\tY\\
&&Q_i^{(d)}\arrow[dd,|->]&~  &\arrow[ll,|->] P_{i,1}^{(d)} \arrow[r,no head,dashed]\arrow[dd,|->] & P_{i,j_i}^{(d)}\arrow[dd,|->] & &\\
& & & & & & \\
X&&P_i &~ &\arrow[ll,|->]  P_{i,1} \arrow[r,no head,dashed] & P_{i,j_i}&&\tX
\end{tikzcd}
\]
The map $Y\to X$ is clearly étale away from the $Q_{i}^{(m)}$. Denoting by $\TC$ the tangent cones, the construction of $Y$ ensures that there is an isomorphism \[ (\TC_{P_{i,1}^{(m)}}(\tY)\sqcup \TC_{P_{i,2}^{(m)}}(\tY) \sqcup\dots \sqcup \TC_{P_{i,j_i}^{(m)}}(\tY))_{P_{i,1}^{(m)}+\dots+P_{i,j_i}^{(m)}}\xrightarrow{\sim}\TC_{Q_i^{(m)}}(Y).\]
Moreover, since $\tY\to \tX$ is étale, the maps $\TC_{P_{i,j}^{(m)}}(\tY)\to \TC_{P_{i,j}}(\tX)$ are isomorphisms. The following commutative diagram now shows that $\TC_{Q_i^{(m)}}(Y)\to \TC_{P_i}(X)$ is an isomorphism for any $i\in \{1\dots r\}$ and $m\in \{1\dots j_i\}$, i.e. that $Y\to X$ is étale.  

\[
\begin{tikzcd}
(\TC_{P_{i,1}^{(m)}}(\tY)\sqcup \TC_{P_{i,2}^{(m)}}(\tY) \sqcup\dots \sqcup \TC_{P_{i,j_i}^{(m)}}(\tY))_{P_{i,1}^{(m)}+\dots+P_{i,j_i}^{(m)}} 
\arrow[r,"\sim"]\arrow[d,"\sim"{anchor=north, rotate=90}] & \TC_{Q_i^{(m)}}(Y)\arrow[d] \\
(\TC_{P_{i,1}}(\tX)\sqcup \TC_{P_{i,2}}(\tX) \sqcup\dots \sqcup 
\TC_{P_{i,j_i}}(\tY))_{P_{i,1}+\dots+P_{i,j_i}} \arrow[r,"\sim"] & \TC_{P_i} (X)
\end{tikzcd}
\]
Furthermore, since the equivalence classes of the quotient map $\tY\to Y$ are stable under the action of $\Aut(\tY|\tX)$, any $\tX$-automorphism of $\tY$ induces an $X$-automorphism of $Y$, and $Y\to X$ is Galois. 
\end{proof}
\end{Lemma}

\begin{Remark} Note that the covering $Y$ of $X$ constructed in the above lemma is irreducible, since its normalisation is connected. The lemma shows that normalisation defines an essentially surjective functor between the category of irreducible Galois coverings of $X$ with automorphism group $G$ and the category of Galois coverings of $\tX$ with automorphism group $G$. 
\end{Remark}

\subsection{Coverings of reducible curves with multicross singularities}

Consider a connected reduced curve $X$ over an algebraically closed field $k$. Suppose that $X$ has at worst multicross singularities. Denote by $X_1,\dots,X_r$ its irreducible components. Since its singularities are multicross, the curve $X$ may be defined by simple divisors $m_1,\dots,m_N\in \Div(\bigsqcup_j X_j)$ with disjoint supports. We may assume that for each $i\in \{1\dots N\}$, there is a subset $J_i$ of $\{ 1\dots r\}$ and for any $j\in J_i$, one point $Q_{i,j}\in X_j$ such that  \[m_i=\sum_{j\in J_i}Q_{i,j}.\]  

\begin{Lemma}\label{lem:redcover} Consider an irreducible component $X_c$ of $X$, and a Galois covering $X_c'\to X_c$ with group $G=\{ \sigma_1=\id,\sigma_2,\dots,\sigma_d\}$. There is a Galois covering $X'\to X$ with an $X_c$-isomorphism $X'\times_X X_c \xrightarrow{\sim}X'_c$.
\begin{proof} For each $j\in \{1\dots r\}- \{c\}$, consider $d$ copies $X_j^1,\dots,X_j^d$ of $X_j$.  We are going to define $X'$ by moduli on $Y\coloneqq X'_c\sqcup\bigsqcup_{j\neq i}\bigsqcup_{s=1}^d X_j^{(s)}$. Pick $i\in\{1\dots N\}$. If $c\in J_i$, choose a preimage $R_{i,c}$ of $Q_{i,c}$ in $X_c'$. For each $s\in \{1\dots d\}$ and $j\in J_i$, set 
\[ Q_{i,j}^{(s)}\coloneqq\left\lbrace \begin{array}{lcl} \sigma_s(R_{i,c}) & \text{if}\qquad & j=c \\
\text{the copy of }Q_{i,j}\text{ lying on } X_j^{(s)}\qquad\qquad & \text{if} \qquad & j\neq c\end{array}\right.\]
For each $s\in\{1\dots d\}$, set \[ m_i^{(s)}\coloneqq \sum_{j\in J_i} Q_{i,j}^{(s)}.\]
Let $X'$ be the curve defined by the moduli $m_1^{(1)},\dots,m_1^{(d)},\dots,m_N^{(1)},\dots,m_N^{(d)}\in \Div(Y)$. Denote by $P_i^{(s)}$ (resp. $P_i$) the image of $m_i^{(s)}$ (resp. $m_i$) in $X'$ (resp. $X$). The following diagram sums up the notations.

\[
\begin{tikzcd}
X'&& P_i^{(1)} \arrow[d,no head,dashed] &~ &\arrow[ll,|->] Q_{i,1}^{(1)}\arrow[d,no head,dashed] \arrow[r,no head,dashed] & Q_{i,|J_i|}^{(1)}\arrow[d,no head,dashed]&&Y\\
&&P_i^{(d)}\arrow[dd,|->]&~  &\arrow[ll,|->] Q_{i,1}^{(d)} \arrow[r,no head,dashed]\arrow[dd,|->] & Q_{i,|J_i|}^{(d)}\arrow[dd,|->] & &\\
& & & & & & \\
X&&P_i &~ &\arrow[ll,|->]  Q_{i,1} \arrow[r,no head,dashed] & Q_{i,|J_i|}&&\bigsqcup_{j\in J_i} X_j
\end{tikzcd}
\]
The curve $X'$ is connected because $X$ is connected: any of the $X_j$ contains a preimage of a singular point of $X$, hence any of the $X_j^{(s)}$ contains a point which lies in the support of one of the $m_i^{(s)}$, whose image $P_i^{(s)}$ lies in the connected component of $X_c'$. The construction of $X'$ ensures that the maps $X_i'\to X_i$ and $X_j^{(s)}\xrightarrow{\sim}X_j$ define a map $X'\to X$. There are $X$-automorphisms $\Sigma_1,\dots,\Sigma_d$ of $X'$ defined, for $s\in \{1\dots d\}$, by $\Sigma_s|_{X_c}=\sigma_s$ and $\Sigma_s|_{X_j^{(t)}}\colon X_j^{(t)}\xrightarrow{\sim} X_j^{(r)}$ where $r$ is defined by $\sigma_s\sigma_t=\sigma_r$. Hence $X'\to X$ is Galois. Finally, let us prove that $X'\to X$ is étale. It is clearly étale away from the intersection points of the irreducible components. Take $i\in \{1\dots N \}$ and $s\in \{1\dots d\}$. 
Consider the following commutative diagram, where we also denote $X_c'$ by $X_c^{(s)}$:
\[
\begin{tikzcd}
\left(\bigsqcup_{j\in J_i}\TC_{Q_{i,j}^{(s)}}(X_j^{(s)})\right)_{\sum_jQ_{i,j}^{(s)}}\arrow[r,"\sim"]\arrow[d] & \TC_{P_i^{(s)}}(X')\arrow[d] \\
\left(\bigsqcup_{j\in J_i}\TC_{Q_{i,j}}(X_j)\right)_{\sum_jQ_{i,j}}\arrow[r,"\sim"]&\TC_{P_i}(X)
\end{tikzcd}
\]
Since the downwards map on the left hand side is an isomorphism by étaleness of $X_c'\to X_c$ and $X_j^{(s)}\xrightarrow{\sim}X_j$, the map $X'\to X$ is étale at $P_i^{(s)}$.
\end{proof}
\end{Lemma}

\subsection{A covering which trivialises $H^1$ and $H^2$ over algebraically closed fields}

Let $n$ be a positive integer. We will denote by $\Lambda$ the ring $\ZZ/n\ZZ$. Given a $\Lambda$-module $M$, we will denote by $M^\vee$ its $\Lambda$-dual. We will only work with profinite groups, so every module we consider will be endowed with the discrete topology.\\

Recall that given a connected noetherian scheme $X$ such that $\HH^1(X,\Lambda)$ is finite, there is a Galois covering $\Xnt$ of $X$ with automorphism group $\HH^1(X,\Lambda)^\vee$. This covering is unique up to isomorphism and corresponds to the quotient of the étale fundamental group $\pi$ of $X$ at some geometric point by the closure of $\pi^n[\pi,\pi]$. For more details, see \cite[Cor. 2.2]{monarticle}. In case $n=\ell$ is a prime number, this corresponds to the Frattini quotient of the pro-$\ell$ completion of $\pi$, which is used in \cite[§3]{mo} to establish the computability of some affine curve fibrations.

\begin{Lemma}\label{lem:morphH1zero} Let $X$ be a connected noetherian scheme such that $\HH^1(X,\Lambda)$ is finite. Then the map \[ \HH^1(X,M)\to \HH^1(X\nt,M)\] 
is trivial.
\begin{proof}
Using the structure of finite abelian groups, it suffices to show this for $M=\Lambda$. Let $\pi$ denote the étale fundamental group of $X$ at some geometric point. Elements of $\HH^1(X,\Lambda)$ are continuous homomorphisms $\pi\to \Lambda$, which are necessarily trivial on $\pi^n[\pi,\pi]$, hence \[\Hom_{\rm cont}(\pi,\Lambda)=\HH^1(X,\Lambda)\to \HH^1(X\nt,\Lambda)\] is trivial. 
\end{proof}
\end{Lemma}

\begin{prop}\label{prop:morphH1zero} Let $X$ be a connected reduced algebraic curve over an algebraically closed field whose characteristic does not divide $n$. Let $M$ be a $\Lambda$-module. If the geometric genus of every irreducible component of $X$ is nonzero, the map \[ \HH^2(X,M)\to \HH^2(X\nt,M)\]
is trivial.
\begin{proof} 

Denote by $X_1,\dots,X_r$ the irreducible components of $X$. Suppose that the geometric genus of each $X_i$ is nonzero. Denote by $\nu\colon \tX\to X$ the normalisation of $X$. Its connected components are the normalisations $\tX_1,\dots,\tX_r$ of the irreducible components of $X$. 
Since the cokernel of $\Lambda\to \nu_\star \nu^\star \Lambda$ is supported on a zero-dimensional subscheme of $X$, the map \[ \HH^2(X,\Lambda)\to \HH^2(X,\nu_\star\nu^\star \Lambda)\simeq\HH^2(\tX,\Lambda)\simeq\prod_{i=1}^r\HH^2(\tX_i,\Lambda)\] is an isomorphism. Let $i\in \{1\dots r\}$. 
Since $\tX_i$ has nonzero genus, the degree of the map $\tX_i\nt\to\tX_i$ is divisible by $n$, and the map $\HH^2(\tX_i,\Lambda)\to\HH^2(\tX_i\nt,\Lambda)$ is multiplication by this degree: it is the trivial map. \Cref{lem:seminormcover} ensures there is a Galois covering $X_i'\to X_i$ such that $X_i'\times_{X_i}\tX_i=\tX_i\nt$. \Cref{lem:redcover} then ensures there is a Galois covering $Y_i\to X$ with the same automorphism group such that $Y_i\times_X X_i=X_i'$. Since $X_i$ is an irreducible component of $X$, we know that $\HH^2(X_i,\Lambda)$ is a direct factor of $\HH^2(X,\Lambda)$. Likewise, $\HH^2(X'_i,\Lambda)$ is a direct factor of $\HH^2(Y_i,\Lambda)$, and the map $\HH^2(X,\Lambda)\to \HH^2(Y_i,\Lambda)$ restricts to $\HH^2(X_i,\Lambda)\to \HH^2(X'_i,\Lambda)$. The following commutative diagram
\[ 
\begin{tikzcd}
\HH^2(\tX_i,\Lambda)\arrow[d,"0"] &\arrow[l,"\sim",swap] \HH^2(X_i,\Lambda)\arrow[r,hookrightarrow]\arrow[d] & \HH^2(X,\Lambda) \arrow[d] \\
\HH^2(\tX_i\nt,\Lambda) &\arrow[l,"\sim",swap] \HH^2(X'_i,\Lambda)\arrow[r,hookrightarrow] & \HH^2(Y_i,\Lambda)
\end{tikzcd}
\]
shows that the restriction to $\HH^2(X_i,\Lambda)$ of the map \[ \HH^2(X,\Lambda)\to\HH^2(Y_i,\Lambda)\]
is trivial.
The group $\HH^1(\tX_i,\Lambda)^\vee$ is a quotient of $\HH^1(\tX,\Lambda)^\vee$, which itself is a quotient of $\HH^1(X,\Lambda)^\vee$. The latter quotient map may be taken to be the dual of a section of the surjective map $\HH^1(X,\Lambda)\to\HH^1(\tX,\Lambda)$ of $\Lambda$-modules, whose target is free. Therefore, since $Y_i\to X$ has automorphism group $\HH^1(\tX_i,\Lambda)^\vee$, the covering $X\nt\to X$ factors through $Y_i\to X$, and trivialises the image of $\HH^2(X_i,\Lambda)$ in $\HH^2(X,\Lambda)$. Since $\HH^2(X,\Lambda)$ is the direct product of these images, this concludes the proof. 
\end{proof}
\end{prop}

\begin{Definition}
Let $\Lambda$ be a ring. A sheaf of $\Lambda$-modules on the (small) étale site of a scheme $X$ will be called \textit{lisse} if it is locally constant and constructible. We say that a connected étale covering $Y\to X$ trivialises such a sheaf $\F$ if the pullback $\F|_Y$ is a constant sheaf.
\end{Definition}

\begin{Corollary}\label{cor:morphH1zero} Let $X$ be a connected reduced algebraic curve with at worst multicross singularities over an algebraically closed field whose characteristic does not divide $n$. Let $\LL$ be a lisse sheaf of $\Lambda$-modules on $X$, and $Y\to X$ be a Galois covering which trivialises $\LL$. Then the map \[ \HH^1(X,\LL)\to \HH^1(Y\nt,\LL|_{Y\nt})\] 
is trivial, and if the geometric genus of every irreducible component of $X$ is nonzero, the map \[ \HH^2(X,\LL)\to \HH^2(Y\nt,\LL|_{Y\nt})\]
is trivial as well.
\begin{proof} The covering $Y\to X$ being étale, $Y$ still only has multicross singularities. Moreover, if the irreducible components of $X$ have nonzero genus, those of $Y$ verify the same property.  Therefore, we may apply \Cref{lem:morphH1zero} and \Cref{prop:morphH1zero} to the sheaf $\LL|_Y$, and for $i=1,2$, the composite map
\[ \HH^i(X,\LL)\to \HH^i(Y,\LL|_Y)\to \HH^i(Y\nt,\LL|_{Y\nt})\]
is trivial.
\end{proof}
\end{Corollary}

\subsection{Explicit construction of this covering for integral curves}\label{sec:expl}

Let $X$ be a projective integral curve with multicross singularities over an algebraically closed field $k$. Let $n$ be an integer invertible in $k$. We denote by $\Lambda$ the ring $\ZZ/n\ZZ$. The covering $X\nt\to X$ described in the previous section is built from two types of coverings: \begin{itemize}[label=$\diamond$]
\item the irreducible covering $\tX\nt$ of the normalisation $\tilde X$ of $X$;
\item reducible coverings built by gluing copies of $\tX$ along preimages of the singular points.
\end{itemize}
The covering $\tX\nt$ is classically constructed by adding to the function field of $\tX$ the $n$-th roots of functions whose divisor is a multiple of $n$; for an in-depth study of its construction and properties, see \cite[§2]{monarticle}.
Let us now turn to the coverings coming from the singularities of $X$. This is a generalisation of the well-known fact that the nodal cubic has a Galois covering with group $\Lambda$ constructed by gluing $n$ copies of $\PP^1$ in a cyclic way along the two preimages of the singular point, as described in \cite[III, Exercise 10.6]{hartshorne}. Let us first introduce some notations. \begin{itemize}[label=$\diamond$]
\item Let $P^{(1)},\dots,P^{(r)}$ be the singular points of $X$. In the next points of this list, we fix $i\in \{1\dots r\}$.
\item We denote by $Q_{0}^{(i)},\dots,Q_{m_i-1}^{(i)}$ the preimages of $P^{(i)}$ in $\tX$, and by $\tY_i$ the scheme $X\times\Lambda^{m_i-1}$. 
\item Given $j\in\{0\dots m_i-1\}$ and $s=(s_1,\dots,s_{m_i-1})\in\Lambda^{m_i-1}$, we denote by $Q^{(i)}_{j,s}$ the point $Q^{(i)}_j$ lying on the connected component of index $s$ in $\tY_i$.
\item For each $s\in \Lambda^{m_i-1}$, consider the divisor \[ D_{s,i}^{(i)}=Q^{(i)}_{0,s}+\sum_{j=1}^{m_i-1}Q^{(i)}_{j,(s_1,\dots,s_{j-1},s_j+1,s_{j+1},\dots,s_{m_i-1})}\in\Div(\tY_i). \]
as well as, for $t\in \{1\dots r\}-\{ i\}$ :
\[ D_{s,i}^{(t)}=\sum_{j=0}^{m_t-1}Q^{(t)}_{j,s}\in\Div(\tY_i).\]
\item Let $Y_i$ be the scheme obtained from $\tY_i$ by contracting each of the divisors $D^{(t)}_{s,i}$ for $s\in \Lambda^{m_i-1}$ and $t\in \{1\dots r\}$ as described in \Cref{subsec:notations}. It comes equipped with a map $Y_i\to X$ sending the image of $D_{s,i}^{(t)}$ to $P^{(t)}$. 
\end{itemize}

\begin{figure}[h]
\begin{center}\includegraphics[scale=0.4]{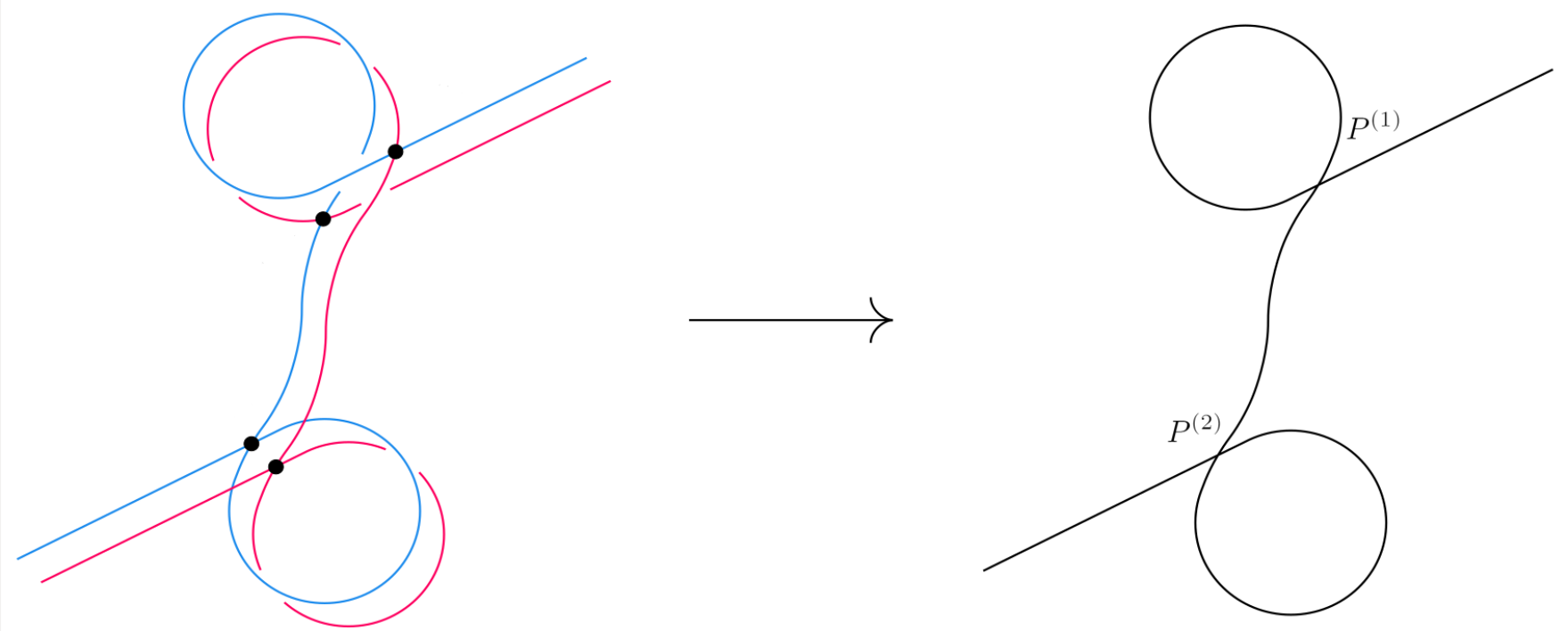}\end{center}
\caption{The covering $Y_1\to X$ when $r=2, m_1=m_2=2$}
\end{figure}

\begin{Lemma} For each $i\in \{1\dots r\}$, the map $Y_i\to X$ is étale with automorphism group $\Lambda^{m_i-1}$.
\begin{proof} Consider $t\in \{1\dots r\}$. The map of tangent cones $\TC_{D_{s,t}^{(i)}}Y_i\to \TC_{P^{(i)}}X$ is nothing more than \[\bigsqcup_{j=0}^{m_i-1}\TC_{Q_{j}^{(i)}}\tX \to \TC_{P^{(i)}}X\]
which itself is an isomorphism since $P^{(i)}$ is a multicross singularity.
Finally, each $\sigma\in\Lambda^{m_i-1}$ defines an $X$-automorphism of $\tX\times\Lambda^{m_i-1}$ by permuting the factors of the product. This automorphism fixes the divisors $D_{s,t}^{(i)}$, hence defines an $X$-automorphism of $Y_i$.
\end{proof}
\end{Lemma}

Recall that $\HH^1(X,\mu_n)$ is canonically isomorphic to the $n$-torsion of the Picard group of $X$. There is a short exact sequence \cite[V, 13, Th. 3]{serrega}:
\[
0\to \bigoplus_{i=1}^r \frac{(k^\times)^{m_i}}{k^\times}\to \Pic^0(X)\to\Pic^0(\tX)\to 0
\]
Passing to the $n$-torsion, one gets: \[ 0\to \frac{\mu_n(k)^{m_1}}{\mu_n(k)}\times \dots\times \frac{\mu_n(k)^{m_r}}{\mu_n(k)} \to \HH^1(X,\mu_n)\to \HH^1(\tX,\mu_n)\to 0.\]
In particular, since $\HH^1(\tX,\mu_n)$ is free, $\HH^1(X,\Lambda)$ is a free $\Lambda$-module isomorphic to $\Lambda^{m_1-1}\times\cdots\times\Lambda^{m_r-1}\times\HH^1(X,\Lambda)$.\\

Consider the scheme $Y\coloneqq Y_1\times_X\dots\times_X Y_r$ as well as the covering $W$ of $X$ with group $\HH^1(\tX,\Lambda)^\vee$ obtained from $\tX\nt\to \tX$ as described in \Cref{lem:seminormcover}, and set $X'=W\times_X Y$.

\begin{prop} The covering $X'$ of $X$ is isomorphic to $\Xnt\to X$.
\begin{proof} It suffices to show that it is a Galois covering of $X$ with Galois group $\HH^1(X,\Lambda)^\vee\simeq \HH^1(\tX,\Lambda)^\vee\times L$, where $L=\Lambda^{m_1-1}\times\dots\times\Lambda^{m_r-1}$. 
Let us first prove that $X'$ is connected. Since $X'\to Y$ is smooth, the normalisation of $X'$ is \[X'\times_Y \tY\simeq W\times_X Y\times_Y \tY\simeq W\times_X \tX\times L\simeq \tX\nt\times L.\]
Consider $s=(s_{i,j})_{i\in \{1\dots r\},j\in \{1\dots m_i-1\}}\in L$. For one such index $(i,j)$, denote by $s'$ the tuple obtained from $s$ by adding 1 to $s_{i,j}$. Since $\Lambda$ is cyclic, to prove that $X'$ is connected, it is enough to show that the intersection of the images in $X'$ of the copies of $\tX\nt$ indexed by $s$ and $s'$ is nonempty. Consider a point $R$ of $W$ above $P^{(i)}$, and the point $R'=(D_{s,1}^{(i)},\dots,D_{s,r}^{(i)})$ of $Y=Y_1\times_X\times\dots\times_XY_r$. Then the point $(R,R')\in X'=W\times_X Y$ lies in said intersection.
Let us now turn to the automorphism group of $X'\to X$. Since each $Y_i\to X$ is surjective, the map $\Aut(Y_1|X)\times\dots\times \Aut(Y_r|X)\to \Aut(Y|X)$ is injective. The same reasoning shows that $\Aut(Y|X)\times\Aut(W|X)\to \Aut(X'|X)$ is injective. Hence $\Aut(X'|X)$ contains $L\times \HH^1(\tX,\Lambda)^\vee$. Since the degree of $X'\to X$ is precisely $|L|\cdot |\HH^1(\tX,\Lambda)|$, this shows that $X'\to X$ is indeed a Galois covering; étaleness follows from that of $Y\to X$ and $W\to X$.
\end{proof}
\end{prop}

\subsection{Trivialising $H^1$ and $H^2$ over fields of cohomological dimension 1}\label{subsec:trivfin}

Let $k$ be a perfect field, and $\bar k$ be an algebraic closure of $k$. Let $n$ be an integer invertible in $k$, and $\Lambda$ be the ring $\ZZ/n\ZZ$. We require $\HH^1(k,\Lambda)$ to be finite, which is true for instance when $k$ is finite, in which case the covering $\Spec(k)\nt$ is the spectrum of the degree $n$ extension of $k$ inside $\bar k$. \\

Let $X$ be a reduced geometrically connected curve over $k$, such that every irreducible component of $X\times_k\bar k$ has nonzero genus. Let $Y$ be a Galois covering of $X$. In this section, we describe a Galois covering $Y\ns$ of $X$ which trivialises $\HH^1$ and $\HH^2$ of constant sheaves on $Y$.\\

We construct the covering $Y\ns$ of $X$ as follows. Let $\bar C$ be a connected component of $\bar Y\coloneqq Y\times_k\bar k$, and $\bar T=\bar C\nt$ be its covering with group $\HH^1(\bar C,\Lambda)^\vee$. Let $K$ be the field of definition of $\bar T$, and $T_K$ be a covering of $X$ defined over $K$ such that $T_K\times_K \bar k$ is isomorphic to $\bar T$. Let $K'$ be the field of definition of $\HH^1(\bar T,\Lambda)$, i.e. the smallest extension of $K$ such that the action of $\Gal(\bar k|K)$ on $\HH^1(\bar T,\Lambda)$ factors through $\Gal(K'|K)$. Denote by $L$ the extension ${K'}\nt$ of $K'$ with group $\HH^1(K',\Lambda)^\vee$. We define $Y\ns\to X$ as the Galois closure of $T_L\coloneqq T_K\times_K L$ over $X$, and denote by $\overline{Y\ns}$ its base change to $\bar k$.

\begin{Remark}\begin{enumerate}
\item Note that $T_K\to X$ factors through $Y_K$. While $T_K\to Y_K$ is surjective only when $Y$ is geometrically connected, the composite map $T_K\to Y_K\to Y$ is always a covering, and $\overline{Y\ns}\to Y$ is a Galois covering.
\item The construction of $T$ only depends on $Y$, but that of $Y\ns$ depends on $X$, because we require $Y\ns$ to be Galois over $X$. 
\end{enumerate}
\end{Remark}

\begin{Lemma}\label{lem:H1H2fin} Let $M$ be the constant sheaf on $Y$ associated to a finite $\Lambda$-module. If $k$ has cohomological dimension 1 then for all $i\in \{1,2\}$ the maps $\HH^i(Y,M)\to \HH^i(Y\ns,M)$ and $\HH^i(\bar Y,M)\to \HH^i(\overline{Y\ns},M)$ are trivial.
\begin{proof} Concerning the maps $\HH^i(\bar Y,M)\to \HH^i(\overline{Y\ns},M)$, it suffices to note that $\overline{Y\ns}\to Y$ factors through $\bar T=\bar C\nt\to \bar C$, which trivialises $\HH^1$ and $\HH^2$ of any constant sheaf by \Cref{cor:morphH1zero}. Let us now turn our attention to $\HH^i(Y,M)\to \HH^i(Y\ns,M)$.
The action of $\Gal(\bar k|K)$ on $\HH^0(\bar T,M)$ is trivial, and $L$ contains $K\nt$, so by \Cref{lem:morphH1zero} $\HH^1(K,\HH^0(\bar T,M))\to \HH^1(L,\HH^0(\bar T,M))$ is trivial as well. The following commutative diagram with exact rows now shows that $\HH^1(Y,M)\to \HH^1(Y\ns,M)$ is trivial.
\[ \begin{tikzcd} 
 ~& \HH^1(Y,M) \arrow[d]\arrow[r]& \HH^0(k,\HH^1(\bar Y,M))\arrow[d,"0"] \\
\HH^1(K,\HH^0(\bar T,M)) \arrow[d,"0",swap]\arrow[r]& \HH^1(T_K,M) \arrow[d]\arrow[r]& \HH^0(K,\HH^1(\bar T,M)) \\
 \HH^1(L,\HH^0(\bar T,M)) \arrow[r]& \HH^1(T_L,M)\arrow[d] & ~ \\
 ~ & \HH^1(Y\ns,M) & ~
\end{tikzcd}
\]
Since $k$ has cohomological dimension 1, so do $K$ and $L$. The Hochschild-Serre spectral sequence \[ E_2^{pq}=\HH^p(K,\HH^q(\bar T_K,M))\Rightarrow \HH^{p+q}(T_K,M) \] degenerates on the second page, and the map \[\ker(\HH^2(T_K,M)\to\HH^0(K,\HH^2(\bar T,M)))\to \HH^1(K,\HH^1(\bar T,M))\] is an isomorphism. This also applies to the Hochschild-Serre spectral sequence of the pullback of $M$ to $T_L$. We obtain the following commutative diagram with exact rows.
\[ \begin{tikzcd} 
 ~& \HH^2(Y,M) \arrow[d]\arrow[r]& \HH^0(k,\HH^2(\bar Y,M))\arrow[d,color=purple] \\
 \HH^1(K,\HH^1(\bar T,M)) \arrow[d,color=cyan]\arrow[r]& \HH^2(T_K,M) \arrow[d]\arrow[r]& \HH^0(K,\HH^2(\bar T,M)) \\
 \HH^1(L,\HH^1(\bar T,M)) \arrow[r]& \HH^2(T_L,M)\arrow[d] & ~ \\
 ~ & \HH^2(Y\ns,M) & ~
\end{tikzcd}
\]
The map $\HH^2(\bar Y,M){\color{purple}\to} \HH^2(\bar T,M)$ factors through $\HH^2(\bar C,M)\to \HH^2(\bar T,M)$, which is trivial since $\bar T=\bar C\nt$. Moreover, by construction of $L$, the map $\HH^1(K,\HH^1(\bar T,M)){\color{cyan}\to} \HH^1(L,\HH^1(\bar T,M))$ is trivial as well, which concludes.
\end{proof}
\end{Lemma}

\begin{Remark}\label{rk:bigcov} The same method may be applied again to the covering $Y\ns\to X$ in order to construct a Galois covering $Y\nns$ of $X$ such that $\HH^i(Y\ns,M)\to\HH^i(Y\nns,M)$ is trivial for $i=1,2$. In that manner, we can construct recursively a tower of Galois coverings \[ Y^{[n]\dots [n]} \to \dots \to Y\ns \to X \]
which we will use in \Cref{subsec:cupfin} to compute cup products. Denoting by $L_i$ the field of definition of the $i$-th cover constructed in this way, the construction above ensures that $L_{i+1}$ contains $L_i\nt$.
\end{Remark}

\section{Most curves are $K(\pi,1)$ schemes}\label{sec:kpi1}

\subsection{Definitions}

We first recall the Definition of a $K(\pi,1)$ scheme. A detailed study of their properties can be found for instance in \cite{achinger_phd}. Let $X$ be a connected noetherian scheme, and $\bar{x}$ a geometric point of $X$. Let $\ell$ be a prime number invertible on $X$. Denote by $\Xet$ the (small) étale topos on $X$. Let $\Xfet$ be the finite étale topos on $X$, i.e. the topos on the site defined by the category of finite étale $X$-schemes endowed with the étale topology. The usual equivalence of categories between finite étale $X$-schemes and continuous finite $\pi_1(X,\bar x)$-sets \cite[VI, §9.8]{abbes_gros} defines an equivalence of topoi between $\Xfet$ and the topos $B_{\pi_1(X,\bar x)}$ of continuous finite $\pi_1(X,\bar x)$-sets.

\begin{Definition} A sheaf $\F$ on $X$ is said to be $\ell$-monodromic if it is lisse and if the image of $\pi_1(X,\bar x)$ in $\Aut(\F_{\bar x})$ is an $\ell$-group. A sheaf of abelian groups $\F$ on $X$ is said to be $\ell$-monodromic if its fibres are finite $\ell$-groups and if it is $\ell$-monodromic as a sheaf of sets.
\end{Definition}

Denote by $\Xlet$ the topos on the site defined by the full subcategory of finite étale schemes $Y\to X$ such that the sheaf on $X$ represented by $Y$ is $\ell$-monodromic. Consider the morphisms of topoi $\rho\colon\Xet\to\Xfet$ and $\rho_\ell\colon\Xet\to\Xlet$ defined by the restriction of sheaves to the sites defining $\Xfet$ et $\Xlet$. Finally, denote by $\pi_1(X,\bar x)^\ell$ the pro-$\ell$ completion of $\pi_1(X,\bar x)$.

\begin{Definition}\cite[Def. 9.21]{abbes_gros} The scheme $X$ is called a $K(\pi,1)$ scheme if for any integer $n$ invertible on $X$ and any sheaf $\F$ of $\ZZ/n\ZZ$-modules on $\Xfet$, the adjunction map \[ \F \to \R\rho_\star(\rho^\star \F)\]
is an isomorphism. It is called a pro-$\ell$ $K(\pi,1)$ scheme if for any $\ell$-primary torsion abelian sheaf $\F$ on $\Xlet$, the adjunction morphism \[ \F\to \R{\rho_\ell}_\star(\rho_\ell^\star\F) \]
is an isomorphism.
\end{Definition}

\begin{Lemma}\cite[§1.4.2]{mo}
The scheme $X$ is a $K(\pi,1)$ scheme if and only if, for any integer $n$ invertible on $X$ and any lisse sheaf $\LL$ of $\ZZ/n\ZZ$-modules, the morphism 
\[\RG(\pi_1(X,\bar x),\LL_{\bar x})\to\RG(X,\LL)\]
is an isomorphism. The scheme $X$ is a pro-$\ell$ $K(\pi,1)$ scheme  if and only if, for any abelian $\ell$-monodromic sheaf $\LL$ on $X$, the morphism \[ \RG(\pi_1(X,\bar x)^\ell,\LL_{\bar x})\to \RG(X,\LL)\]
is an isomorphism.
\end{Lemma}

The fact that $X$ is a $K(\pi,1)$ scheme is equivalent to the composition $\RG(X,\rho^\star -)$ being the derived functor of $\Gamma(\pi_1(X,\bar x),-)$. The following proposition restates this in terms of effaceability.

\begin{prop}\cite[Prop. 2.1.5]{achinger_phd}\label{prop:condKpi1} Let $X$ be a connected scheme and $\ell$ be a prime number invertible on $X$. The scheme $X$ is a $K(\pi,1)$ (resp. pro-$\ell$ $K(\pi,1)$) scheme if and only if for every integer $i>0$ and every lisse sheaf $\LL$ with torsion invertible on $X$ (resp. lisse sheaf of finite $\ZZ/\ell\ZZ$-vector spaces), there exists a Galois covering $\phi_i\colon Y_i\to X$ (with automorphism group an $\ell$-group) such that the map $\HH^i(X,\LL)\to \HH^i(Y_i,\phi_i^\star\LL)$ is trivial.
\end{prop}

\subsection{Proof of the theorem}

The following result is already well-known in the case of smooth or affine curves; it is proven for instance in \cite[Prop. A.4.1]{stix_thesis} using Poincaré duality for smooth curves. As pointed out to the author by Daniel Litt, a similar strategy could be used at least in the case of irreducible curves. Our proof however builds on the explicit constructions of the previous section, which allows full control over the trivialising coverings we use. The following lemma will be used during several reduction steps in the proof.

\begin{Lemma}\label{lem:univhomeo}
Let $f\colon Y\to X$ be a universal homeomorphism of connected schemes. Then $X$ is a $K(\pi,1)$ (resp. a pro-$\ell$ $K(\pi,1)$) if and only if $Y$ is.
\begin{proof} Consider a geometric point $\bar y$ of $Y$, and its image $\bar x$ in $X$.
The result follows directly from the fact that the universal homeomorphism $f\colon Y\to X$ induces isomorphisms between the étale toposes \cite[03SI]{stacks} and between the fundamental groups \cite[0BTT]{stacks} of $X$ and $Y$: given a lisse sheaf $\LL$ on $X$, we have the following commutative diagram.
\[ 
\begin{tikzcd}
\RG(\pi_1(X,\bar x),\LL_{\bar x}) \arrow[d]\arrow[r,"\sim"] & \RG(\pi_1(Y,\bar y),\LL_{\bar y})\arrow[d] \\
\RG(X,\LL) \arrow[r,"\sim"]& \RG(Y,f^\star\LL)
\end{tikzcd}
\] 
\end{proof}
\end{Lemma}

\begin{prop}\label{th:kpi}
Let $X$ be a connected curve over an algebraically closed field $k$. If the geometric genus of every irreducible component of $X$ is nonzero then $X$ is a $K(\pi,1)$ scheme and, for any prime number $\ell$ invertible in $k$, a pro-$\ell$ $K(\pi,1)$ scheme.
\begin{proof} We may suppose that $X$ is reduced, since the canonical map $X_{\rm red}\to X$ is a universal homeomorphism. Let $X^+\to X$ be the seminormalisation map \cite[0EUT]{stacks}. This map is a universal homeomorphism \cite[0EUS]{stacks}, so we may replace $X$ with $X^+$, which has at worst multicross singularities \cite[§2, Cor. 1, (2)]{davis_seminorm}, and assume in the remainder of this proof that $X$ has at worst multicross singularities.
Let $\LL$ be a lisse sheaf of $\ZZ/n\ZZ$-modules on $X$, with $n$ invertible in $k$.  Consider a Galois covering $Y\to X$ which trivialises $\LL$.  \Cref{cor:morphH1zero} ensures that $\HH^i(X,\LL)\to \HH^i(Y\nt,\LL|_{Y\nt})$ is zero for any integer $i\geqslant 1$, and \Cref{prop:condKpi1} concludes. In case $n=\ell$ is prime and the sheaf $\LL$ is $\ell$-monodromic, we may suppose that $\Aut(Y|X)$ is an $\ell$-group. Then $\Aut(Y\lt|X)$ is still an $\ell$-group, which proves that $X$ is also a pro-$\ell$ $K(\pi,1)$.
\end{proof}
\end{prop}

\begin{theorem}\label{th:kpi1} Let $X$ be a geometrically connected curve over a field $k$. Denote by $\bar k$ an algebraic closure of $k$. If the geometric genus of every irreducible component of $X\times_k \bar k$ is nonzero then $X$ is a $K(\pi,1)$ scheme.
\begin{proof} Denote by $\kpf$ the perfect closure of $k$ inside $\bar k$, and by $\Xpf$ the base change of $X$ to $\kpf$. Since $\Xpf\to X$ is induced by base change to a purely inseparable field extension, it is a universal homeomorphism. Hence by \Cref{lem:univhomeo} we may suppose that $k$ is perfect.
Denote by $\GG$ the group $\Gal(\bar k|k)$, and by $\bar X$ the base change of $X$ to $\bar k$. Let $\bar x$ be a generic geometric point of the curve $\bar X$; we will also denote by $\bar x$ its image in $X$. Let $n$ be an integer invertible in $k$. Let $\LL$ be a lisse sheaf of $\ZZ/n\ZZ$-modules on $X$. 
Being the quotient of $\pi_1(X,\bar x)$ by $\pi_1(\bar X,\bar x)$ \cite[IX, Th. 6.1]{sga1}, the group $\GG$ acts on $\RG(\pi_1(\bar X,\bar x),\LL_{\bar x})$. The map $\RG(\pi_1(\bar X,\bar x),\LL_{\bar x})\to \RG(\bar X,\LL|_{\bar x})$ is $\GG$-equivariant, and we obtain the following commutative diagram: 
\[
\begin{tikzcd}
\RG(\pi_1(X,\bar x),\LL_{\bar x}) \arrow[d]\arrow[r,"\sim"] & \RG(\mathfrak{G},\RG(\pi_1(\bar X,\bar x),\LL_{\bar x}))\arrow[d] \\
\RG(X,\LL) \arrow[r,"\sim"]& \RG(\mathfrak{G},\RG(\bar X,\LL|_{\bar X}))  
\end{tikzcd}
\]
\Cref{th:kpi} shows that the downwards map on the right hand side is an isomorphism, which proves that $X$ is a $K(\pi,1)$ scheme. 
\end{proof} 
\end{theorem}
We now know that the cohomology of a lisse sheaf on a curve as in the theorem is that of the corresponding $\pi_1$-module. In the next two sections, we are going to make this more explicit, and replace this profinite group cohomology with the cohomology of finite groups corresponding to explicitly described Galois coverings of the curve.

\section{Computing $\HH^1$, $\HH^2$ and $\HH^3$ in Galois cohomology}\label{sec:galoiscoh}

In the remainder of this article, we show to what extent \Cref{th:kpi1} can be made effective in order to actually compute étale cohomology groups. In practice, people often work with curves over finite or algebraically closed fields. In both these cases, we will construct Galois coverings of a curve, defined using explicitly described quotients of the fundamental group, that allow to compute the cohomology groups of lisse sheaves on it, as well as the cup products between these cohomology groups. Effective computation of cup products has been seen as a promising way of obtaining multilinear maps that might be useful in cryptography (see \cite{boneh} and \cite{chinburg_bleher}). 
In this section, we state a few technical results in Galois cohomology, which we shall apply to fundamental groups of curves in the next one.

\begin{Lemma}\label{lem:H1coh} Let $\pi$ be a profinite group, and $M$ be a finite $\pi$-module. Consider a normal closed subgroup $U$ of $\pi$ such that the restriction map $\HH^1(\pi,M)\to \HH^1(U,M)$ is trivial. Then the inflation map \[\HH^1(\pi/U,\HH^0(U,M))\to \HH^1(\pi,M)\] is an isomorphism.
\begin{proof}
This follows immediately from the low-degree exact sequence of the Hochschild-Serre spectral sequence:
\[ 0\to \HH^1(\pi/U,\HH^0(U,M))\to \HH^1(\pi,M) \to \HH^0(\pi/U,\HH^1(U,M)).\qedhere\] 
\end{proof}
\end{Lemma}

\begin{prop}\label{prop:H2im} Let $\pi$ be a profinite group, and $M$ be a $\pi$-module. Consider normal closed subgroups $U,U',U''$ of $\pi$ such that $U''\subset U'\subset U$. Denote by $G,G',G''$ the respective quotients $\pi/U, \pi/U',\pi/U''$. Suppose that $U$ acts trivially on $\pi$, and that the restriction maps $\HH^2(\pi,M)\to \HH^2(U,M)$, $\HH^1(U,M)\to \HH^1(U',M)$ and $\HH^1(U',M)\to \HH^1(U'',M)$ are all trivial. Then the restriction of the inflation map $\HH^2(G'',M)\to \HH^2(\pi,M))$ to the image of $\HH^2(G',M)\to\HH^2(G'',M)$ is an isomorphism.
\begin{proof}
Consider the morphisms between the Hochschild-Serre spectral sequences (with values in the $\pi$-module $M$) associated with the rows of the following commutative diagram.
\[ 
\begin{tikzcd}
1 \arrow[r] &U'' \arrow[r] \arrow[d] &\pi \arrow[r]\arrow[d] & G'' \arrow[r]\arrow[d] & 1 \\
1 \arrow[r] &U' \arrow[r]\arrow[d] & \pi \arrow[r]\arrow[d] & G' \arrow[r]\arrow[d] & 1 \\
1 \arrow[r] &U \arrow[r] & \pi \arrow[r] & G \arrow[r]& 1 \\
\end{tikzcd}
\]
They yield the following commutative diagram with exact rows, where the group $\HH^2(\pi,M)$ appears in full since the edge maps leaving it factor through the trivial map $\HH^2(\pi,M)\to \HH^2(U,M)$.
\[\begin{adjustbox}{width=\textwidth}{
\begin{tikzcd}0 \arrow[r] & \HH^1(G,M)\arrow[r] \arrow[d] & \HH^1(\pi,M) \arrow[r] \arrow[d] & \HH^0(G,\HH^1(U,M))\arrow[r] \arrow[d] & \HH^2(G,M)\arrow[r] \arrow[d] & \HH^2(\pi,M)\arrow[r] \arrow[d] & \HH^1(G,\HH^1(U,M))\arrow[d] \\
0 \arrow[r] & \HH^1(G',M)\arrow[r] \arrow[d] & \HH^1(\pi,M) \arrow[r,color=cyan] \arrow[d] & \HH^0(G',\HH^1(U',M))\arrow[r,color=green] \arrow[d,color=green] & \HH^2(G',M)\arrow[r,color=green] \arrow[d,color=green] & \HH^2(\pi,M)\arrow[r,color=purple] \arrow[d,color=green] & \HH^1(G',\HH^1(U',M))\arrow[d] \\
0 \arrow[r] & \HH^1(G'',M)\arrow[r] & \HH^1(\pi,M) \arrow[r,color=cyan] & \HH^0(G'',\HH^1(U'',M))\arrow[r,color=green] &  \HH^2(G'',M)\arrow[r,color=green]  & \HH^2(\pi,M)\arrow[r,color=purple]& \HH^1(G'',\HH^1(U'',M))
\end{tikzcd}
}
\end{adjustbox}
\]
The map $\HH^1(U,M)\to \HH^1(U',M)$ is also trivial. Hence the upper right square of the above diagram reads:
\[ \begin{tikzcd}
\HH^2(\pi,M) \arrow[r]\arrow[d,"\sim"{anchor=north, rotate=90}] & \HH^1(G,\HH^1(U,M)) \arrow[d,"0"] \\
\HH^2(\pi,M) \arrow[r,color=purple] & \HH^1(G',\HH^1(U',M))
\end{tikzcd}
\]
and the map $\HH^2(\pi,M){\color{purple}{\to}} \HH^1(G',\HH^1(U',M))$ is trivial. So is $\HH^2(\pi,M){\color{purple}{\to}} \HH^1(G'',\HH^1(U'',M))$. Since the maps $\HH^1(\pi,M){\color{cyan}{\to}}\HH^1(U',M)$ and $\HH^1(\pi,M){\color{cyan}{\to}} \HH^1(U'',M)$ are zero as well, we end up with the following diagram with exact rows.

\[ 
\begin{tikzcd}
0 \arrow[r] & \HH^0(G',\HH^1(U',M)) \arrow[r,color=green] \arrow[d,"0",color=green] & \HH^2(G',M) \arrow[r,color=green] \arrow[d,color=green] & \HH^2(\pi,M) \arrow[r]\arrow[d,"\sim"{anchor=north, rotate=90},color=green] & 0\\
0 \arrow[r] & \HH^0(G'',\HH^1(U'',M)) \arrow[r,color=green] & \HH^2(G'',M) \arrow[r,color=green] & \HH^2(\pi,M) \arrow[r] & 0
\end{tikzcd} 
\]
The map $\HH^2(G',M)\to \HH^2(G'',M)$ defines a map \[ \HH^2(G',M)/\HH^0(G',\HH^1(U',M)) \to \HH^2(G'',M)\]
which, when followed by $\HH^2(G'',M)\to\HH^2(\pi,M)$, becomes an isomorphism. Hence the restriction of the inflation map $\HH^2(G'',M)\to \HH^2(\pi,M)$ to the image of $\HH^2(G',M)\to \HH^2(G'',M)$ is an isomorphism.
\end{proof}
\end{prop}

Even better than the cohomology groups of $M$, we can compute a representative of the cohomology complex $\RG(\pi,M)$ involving only finite groups. Given a group $\mathfrak{S}$ and a $\mathfrak{S}$-module $A$, denote by $C^i(\mathfrak{S},A)$ the group of maps $\mathfrak{S}^i\to A$. We know that $\RG(\mathfrak{S},A)$ is represented by the complex $(C^i(\mathfrak{S},A))_{i\geqslant 0}$ with the usual coboundary maps.\\

\begin{Corollary}\label{cor:RGamma} Consider the situation of \Cref{prop:H2im}. Suppose that $\pi$ has cohomological dimension 2.  Denote by $\psi_2$ the restricion of the inflation map  $C^2(G',M)\to C^2(G'',M)$ to $\ker(C^2(G',M)\to C^3(G',M))$, and by $\partial^1$ the coboundary map $C^1(G'',M)\to C^2(G'',M)$.
The object $\RG(\pi,M)$ in the derived bounded category of abelian groups is represented by the following complex:
\[ M\to C^1(G'',M) \xrightarrow{\partial^1} \im(\partial^1)+\im(\psi_2)\]
where the sum in the last term is the usal sum of submodules inside $C^2(G'',M)$.
\begin{proof}
Denote this complex by $K$. Recall that $\RG(G'',M)$ is represented by the following complex:
\[ M\to C^1(G'',M)\to C^2(G'',M).\]
Let us prove that the natural morphism \[K\to \RG(G'',M)\to \RG(\pi,M) \]
is a quasi-isomorphism. By construction, $K\to \RG(G'',M)$ induces isomorphisms of cohomology groups in degrees 0 and 1. Since $\HH^1(\pi,M)\to \HH^1(U'',M)$ is trivial, \Cref{lem:H1coh} ensures that the inflation map $\HH^1(G'',M)\to \HH^1(\pi,M)$ is an isomorphism. Moreover, the image of $\HH^2(K)\to\HH^2(G'',M)$ is $(\im(\psi_2)+\im(\partial_1))/\im(\partial_1)$, which is exactly the image $\im(\psi_2)/(\im(\partial^1)\cap\im(\psi_2))$ of $\HH^2(G',M)$ in $\HH^2(G'',M)$: \Cref{prop:H2im} concludes.
\end{proof}
\end{Corollary}

\begin{Remark} The above proof is very specific to cohomological dimension 2. Indeed, the cohomology groups in degree 0 and 1 can always be computed as cohomology groups of quotients of $\pi$, which is generally not the case in degree 2 or higher, where we have to restrict to the image of an inflation map as in the proof. This trick only works once, on the last term of the complex.
\end{Remark}

\begin{Lemma}\label{lem:H3coh} Let $\pi$ be a profinite group, and $M$ be a $\pi$-module. Consider a normal closed subgroup $U$ of $\pi$ such that $U$ has cohomological dimension 2 and $\pi/U$ has cohomological dimension 1. Then, in the associated Hochschild-Serre spectral sequence, there is an edge map \[ \HH^3(\pi,M)\to \HH^1(\pi/U,\HH^2(U,M))\]
which is an isomorphism.
\begin{proof}
The conditions we imposed on $U$ and $\pi/U$ imply that the Hochschild-Serre spectral sequence $\HH^p(\pi/U,\HH^q(U,M))\Rightarrow \HH^{p+q}(\pi,M)$ degenerates on page 2. Thus, the quotients in the associated filtration of $\HH^3(\pi,M)$ are the actual cohomology groups $\HH^p(\pi/U,\HH^q(U,M))$ with $p+q=3$. Since there is only one nonzero quotient, namely $\HH^1(\pi/U,\HH^2(U,M))$, the edge map it receives from $\ker(\HH^3(\pi,M)\to \HH^0(\pi/U,\HH^3(U,M))=\HH^3(\pi,M)$ is an isomorphism. 
\end{proof}
\end{Lemma}

\begin{Lemma}\label{lem:H3cohzero} Let $\pi$ be a profinite group of cohomological dimension 2, and $M$ be a $\pi$-module. Consider normal closed subgroups $U''\subset U'\subset U$ of $\pi$ which act trivially on $M$, and denote by $G,G',G''$ the respective quotients $\pi/U,\pi/U',\pi/U''$. If  the restriction maps $\HH^2(U,M)\to\HH^2(U',M)$ and $\HH^1(U',M)\to\HH^1(U'',M)$ are trivial, then the inflation map $\HH^3(G,M)\to \HH^3(G'',M)$ is trivial.
\begin{proof} 
We consider the morphism of Hochschild-Serre spectral sequences:
\[ E_r^{p,q}=\HH^p(G,\HH^q(U,M))\to E_r^{\prime p,q}=\HH^p(G',\HH^q(U',M))\to E_r^{\prime\prime p,q}=\HH^p(G'',\HH^q(U'',M)).\]
Our goal is to prove that $E_2^{3,0}\to E_2^{\prime\prime 3,0}$ is trivial. Since $\pi$ has cohomological dimension $\leqslant 2$, so do its closed subgroups, and all three spectral sequences have only 3 nonzero rows at the second page, and thus degenerate at the fourth page. The abutments of $E$ and $E'$ have no cohomology in degree 3, thus the groups $E_4^{3,0}$ and $E_4^{\prime 3,0}$ are trivial; they are the cokernels of the maps $E_3^{0,2}\to E_3^{3,0}$ and $E_3^{\prime 0,2}\to E_3^{\prime 3,0}$, which are therefore surjective. Since $E_3^{0,2}$ and $E_3^{\prime 0,2}$ are subgroups of $E_2^{0,2}$ and $E_2^{\prime 0,2}$ respectively, and since $\HH^2(U,M)\to\HH^2(U',M)$ is the zero map, the map $E_3^{0,2}\to E_3^{\prime 0,2}$ is trivial as well, and so is $E_3^{3,0}\to E_3^{\prime 3,0}$.
Moreover, the map $E_2^{\prime 1,1}\to E_2^{\prime\prime 1,1}$ is trivial because it is a restriction of $\HH^1(U',M)\to\HH^1(U'',M)$, which itself is trivial. The maps $E\to E'\to E''$ produce the following commutative diagram, whose middle row is exact:
\[
\begin{tikzcd}
~ & E_2^{3,0}\arrow[r]\arrow[d] & E_3^{3,0}\arrow[d,"0"] \\
E_2^{\prime 1,1}\arrow[r]\arrow[d,"0",swap] & E_2^{\prime 3,0}\arrow[d]\arrow[r] & E_3^{\prime 3,0}  \\
E_2^{\prime\prime 1,1}\arrow[r] & E_2^{\prime\prime 3,0} & ~
\end{tikzcd}
\] 
from which we conclude that $E_2^{3,0}\to E_2^{\prime\prime 3,0}$ is the zero map.
\end{proof}
\end{Lemma}

\section{Computing cohomology groups and cup products}\label{sec:cup}

In this section, we use the fact that curves are $K(\pi,1)$ spaces in order to express the cohomology groups of lisse sheaves on them, as well as cup products between these groups, in terms of finite group cohomology. Here, the emphasis lies on constructive results: we have explicitly described in \Cref{sec:Galoiscov} the construction of specific (towers of) Galois coverings of curves, whose automorphism groups will be the only ones considered in this section. 

\subsection{In the cohomology of curves over algebraically closed fields}

Consider a connected curve $X$ over an algebraically closed field $k$. Let $\bar x$ be a geometric point of $X$. Suppose that every irreducible component of $X$ has nonzero genus. \Cref{th:kpi} ensures that $X$ is a $K(\pi,1)$ scheme.  Let $n$ be an integer prime to the characteristic of $k$, and $\LL$ be a lisse sheaf of $\ZZ/n\ZZ$-modules on $X$, with geometric fibre $M$ at $\bar x$. We are going to show how to compute the cohomology of $\LL$ using only finite quotients of $\pi_1(X,\bar x)$.

Let $Y$ be a Galois covering of $X$ such that $\LL|_Y$ is constant. 
Denote respectively by $G^{(1)},G^{(2)},G^{(3)}$ the automorphism groups of $Y\nt\to X$, $Y\nnt\to X$ and $Y\nnnt\to X$.
Denote by $N$ the submodule of $C^2(G^{(3)},M)$ defined as the sum of the images of the maps $C^1(G^{(2)},M)\to C^2(G^{(3)},M)$ and $\ker(C^2(G^{(2)},M)\to C^3(G^{(2)},M))\to C^2(G^{(3)},M)$.

\begin{prop} In this situation, the cohomology complex $\RG(X,\LL)$ is represented by 
\[ M\to C^1(G,M)\to N.\]
\begin{proof}
According to \Cref{cor:morphH1zero}, the subgroups $\pi_1(Y\nt)$, $\pi_1(Y\nnt)$, $\pi_1(Y\nnnt)$ of $\pi_1(X)$ (with suitable choices of geometric points) satisfy the requirements of \Cref{prop:H2im}. \Cref{cor:RGamma} concludes.
\end{proof}
\end{prop}

Now consider two lisse sheaves $\LL_1$ and $\LL_2$ of $\ZZ/n\ZZ$-modules on $X$, both trivialised by the Galois covering $Y\to X$. Denote by $M_1$ and $M_2$ their respective geometric fibres, and by $M$ the tensor product $M_1\otimes M_2$.

\begin{prop}\label{th:cupprod} The cup product
\[ \HH^1(X,\LL_1)\otimes \HH^1(X,\LL_2) \to \HH^2(X,\LL_1\otimes \LL_2) \]
is computed by the following composition:
\[ \HH^1(G^{(2)},M_1)\otimes \HH^1(G^{(2)},M_2) \xrightarrow{\cup} \HH^2(G^{(2)},M) \rightarrow \im(\HH^2(G^{(2)},M)\to \HH^2(G^{(3)},M)).\]
\begin{proof} The coverings $Y\nnnt\to Y\nnt\to Y\nt$ of $X$ satisfy the requirements of \Cref{prop:H2im}, which ensures that the maps
\[ \im (\HH^2(G^{(2)},M) \to \HH^2(G^{(3)},M))\to \HH^2(\pi_1(X),M)\to \HH^2(X,\LL_1\otimes \LL_2) \]
are isomorphisms. Moreover, $\HH^1(G^{(2)},M)\to \HH^1(X,\LL_1\otimes \LL_2)$ is an isomorphism by \Cref{lem:H1coh}. The following commutative diagram:
\[ 
\begin{tikzcd}
\HH^1(G^{(2)},M_1)\otimes \HH^1(G^{(2)},M_2) \arrow[d,"\sim"{anchor=south, rotate=90}] \arrow[r] & \HH^2(G^{(2)},M) \arrow[d] \\
\HH^1(G^{(3)},M_1)\otimes \HH^1(G^{(3)},M_2) \arrow[d,"\sim"{anchor=south, rotate=90}] \arrow[r] & \HH^2(G^{(3)},M) \arrow[d] \\
\HH^1(\pi_1(X),M_1) \otimes \HH^1(\pi_1(X),M_2) \arrow[d,"\sim"{anchor=south, rotate=90}] \arrow[r] & \HH^2(\pi_1(X),M)\arrow[d,"\sim"{anchor=south, rotate=90}] \\
\HH^1(X,\LL_1)\otimes \HH^1(X,\LL_2) \arrow[r] & \HH^2(X,\LL_1\otimes \LL_2)
\end{tikzcd}
\]
where the first three lines are cup products in group cohomology,
shows how to compute the bottom map, which is the cup product in étale cohomology. It is the composition of the maps
\[ \HH^1(G^{(2)},M_1)\otimes \HH^1(G^{(2)},M_2) \xrightarrow{\cup} \HH^2(G^{(2)},M) \rightarrow \im(\HH^2(G^{(2)},M)\to \HH^2(G^{(3)},M)).\qedhere\]
\end{proof}
\end{prop}

\subsection{In the cohomology of curves over fields of cohomological dimension 1}\label{subsec:cupfin}

We now turn our attention to curves defined over a perfect field $k$ of cohomological dimension 1 such that $\HH^1(k,\Lambda)$ is finite. For instance, this could be a finite field, or the perfect closure of the field of fractions of a stricly henselian DVR. In particular, when $k$ is finite, the results in this section give another expression for the cup products computed in \cite{chinburg_bleher}, where the case of constant sheaves is treated in terms of pairings between torsion subgroups of Picard groups. 

Denote by $\bar k$ an algebraic closure of $k$. Let $X$ be a connected curve over $k$ satisfying the conditions of \Cref{th:kpi1}: it is a $K(\pi,1)$ scheme. In the following, the fundamental group of $X$ and every geometric fibre will be with respect to a fixed geometric point of $X$. 
Let $n$ be an integer prime to the characteristic of $k$. As usual, we denote by $\Lambda$ the ring $\ZZ/n\ZZ$. In this section, we describe cup products in the cohomology of locally constant sheaves in terms of the Galois coverings defined in \Cref{subsec:trivfin}.
\\

Consider a lisse sheaf $\LL$ of $\Lambda$-modules on $X$, trivialised by a Galois covering $Y\to X$. Denote by $M$ its geometric fibre. We have described in \Cref{rk:bigcov} a tower of Galois coverings \[ Y_i\ns\coloneqq Y^{{[n]\dots [n]}} \to \dots \to Y_1\ns\coloneqq Y\ns \to X. \] For each integer $i\geqslant 0$, we set $G^{(i)}=\Aut(Y_i\ns|X)$ and define $\bar G^{(i)}$ to be the image of $\pi_1(\bar X)$ in $G^{(i)}$. We denote by $\GG$ the group $\Gal(\bar k|k)$ and by $\SS^{(i)}=\Aut(L_i|k)$ the quotient $G^{(i)}/\bar G^{(i)}$. 
In the following, we will use the homotopy exact sequence \cite[IX, Th. 6.1]{sga1}
\[1 \to \pi_1(\bar X)\to \pi_1(X)\to \GG\to 1 \]
and the quotient sequences 
\[1 \to \bar G^{(i)}\to G^{(i)} \to \SS^{(i)}\to 1. \]
As mentioned in \Cref{rk:bigcov}, $L_{i+1}$ contains the extension $L_i\nt$ of $L_i$ with group $\HH^1(L_i,\Lambda)^\vee$.
According to \Cref{lem:H1H2fin} and \Cref{lem:H3cohzero}, the maps \[ \HH^j(Y_i\ns,M)\to\HH^j(Y_{i+1}\ns,M)\]
and \[ \HH^j(Y_i\ns\times_k \bar k,M)\to \HH^j(Y_{i+1}\ns\times_k \bar k,M)\] are trivial for all $i\geqslant 1$ and $j\in \{1,2\}$.
Thus, by \Cref{lem:H1coh} and \Cref{prop:H2im}:
\begin{itemize}[label=$\diamond$]
\item $\HH^1(G^{(i)},M)\to \HH^1(X,\LL)$ is an isomorphism for all $i\geqslant 1$ ;
\item $\im(\HH^2(G^{(i)},M)\to \HH^2(G^{(j)},M)) \to \HH^2(X,\LL)$ is an isomorphism for all $j>i\geqslant 2$;
\item $\HH^3(G^{(i)},M)\to \HH^3(G^{(j)},M)$ is trivial for all $i,j\geqslant 1$ such that $j\geqslant i+2$.
\end{itemize}
In particular, cup products of the form $\HH^1\times \HH^1\to\HH^2$ may be computed as described in the previous section.

\begin{Lemma}\label{lem:cohhs} Consider integers $i,j\geqslant 1$ such that $j>i\geqslant 2$. The inflation map \[\HH^1(\SS^{(j)},\im(\HH^2(\bar G^{(i)},M)\to\HH^2(\bar G^{(j)},M)))\to \HH^1(\GG,\HH^2(\pi_1(\bar X),M))\]
is an isomorphism. Moreover, there is a canonical isomorphism
\[ \HH^3(X,\LL)\to \HH^1(\GG,\HH^2(\pi_1(\bar X),M)).\]
\begin{proof} \Cref{prop:H2im} ensures that the restriction of $\HH^2(\bar G^{(j)},M)\to\HH^2(\pi_1(\bar X),M)$ to the image of $\HH^2(\bar G^{(i)},M)$ is an isomorphism. The action of $\GG$ on $\im(\HH^2(\bar G^{(i)},M)\to\HH^2(\bar G^{(j)},M)))$ factors through $\SS^{(i)}$, and since the extension $L^{(j)}$ of $L^{(i)}$ contains ${L^{(i)\langle n\rangle}}$, the first statement follows from \Cref{lem:H1coh}. The second statement is obtained by applying \Cref{lem:H3coh} to the subgroup $\pi_1(\bar X)$ of $\pi_1(X)$.
\end{proof}
\end{Lemma}

\begin{prop} Let $\LL_1,\LL_2$ be lisse sheaves on $X$. Denote by $M_1,M_2$ their geometric fibres, and set $M=M_1\otimes M_2$. Let $Y$ be a Galois covering of $X$ such that $\LL_1|_Y$ and $\LL_2|_Y$ are constant. The cup product $\HH^1(X,\LL_1)\otimes \HH^1(X,\LL_2)\to \HH^3(X,\LL_1\otimes \LL_2)$ is computed by the following composition:
\[
\begin{tikzcd}
\HH^1(G^{(3)},M_1)\otimes \im(\HH^2(G^{(2)},M_2)\to \HH^2(G^{(3)},M_2)) \arrow[d,"\cup"]  \\ 
\im(\HH^3(G^{(2)},M)\to\HH^3(G^{(3)},M))\arrow[d,"{\rm infl}"] \\
\ker(\HH^3(G^{(4)},M)\to \HH^3(\bar G^{(4)},M)) \arrow[d,"\delta"] \\ 
\HH^1(\SS^{(4)},\HH^2(G^{(4)},M)) \arrow[d,"{\rm infl}"] \\
\HH^1(\SS^{(5)},\im(\HH^2(G^{(4)},M)\to \HH^2(G^{(5)},M))
\end{tikzcd}
\]
where $\delta$ is an edge map of the Hochschild-Serre spectral sequence associated with the normal subgroup $\bar G^{(4)}$ of $G^{(4)}$.
\begin{proof} Denote by $H^{i,(j)}(-)$ the image of $\HH^i(G^{(j-1)},-)\to \HH^i(G^{(j)},-)$. We know that $H^{1,(3)}(M_1)=\HH^1(G^{(3)},M_1)$ is isomorphic to $\HH^1(X,\LL_1)$, and that $H^{2,(3)}$ is isomorphic to $\HH^2(X,\LL_2)$. The image of $H^{1,(3)}(M_1)\otimes H^{2,(3)}(M_2)$ in $\HH^3(G^{(3)},M)$ lies inside $H^{3,(3)}(M)$. Now by \Cref{lem:H3cohzero} applied to the quotients $\bar G^{(2)},\bar G^{(3)},\bar G^{(4)}$ of $\pi_1(\bar X)$, the map $\HH^3(\bar G^{(2)},M){\color{green}{\to}} \HH^3(\bar G^{(4)},M)$ is trivial. The situation is summed up in the commutative diagram below, from which we conclude that the image of $H^{1,(3)}(M_1)\otimes H^{2,(3)}(M_2){\color{purple}{\to}}H^3(G^{(4)},M)$ lies in the kernel of the restriction map $\HH^3(G^{(4)},M){\color{cyan}{\to}} \HH^3(\bar G^{(4)},M)$. 
\[ 
\begin{tikzcd}
~&\HH^3(G^{(4)},M) \arrow[r,color=cyan] & \HH^3(\bar G^{(4)},M) \\
H^{1,(3)}(M_1)\otimes H^{2,(3)}(M_2) \arrow[r]\arrow[ru,color=purple] &\HH^3(G^{(3)},M)\arrow[u]  & ~ \\
\HH^1(G^{(2)},M_1) \otimes \HH^2 (G^{(2)},M_2) \arrow[r] \arrow[u, two heads] &\HH^3(G^{(2)},M) \arrow[u]\arrow[r] & \HH^3(\bar G^{(2)},M) \arrow[uu,"0",swap,color=green]
\end{tikzcd}
\]
This kernel is the source of an edge map with target $\HH^1(\SS^{(4)},\HH^2(\bar G^{(4)},M))$. The inflation map $\HH^1(\SS^{(4)},\HH^2(\bar G^{(4)},M))\to \HH^1(\SS^{(5)},\HH^2(\bar G^{(5)},M))$ factors through $\HH^1(\SS^{(5)},H^{2,(5)}(M))$. \Cref{lem:cohhs} ensures that this last group is isomorphic to $\HH^1(\GG,\HH^2(\bar X,\LL))$ via an inflation map, which itself is isomorphic to $\HH^3(X,\LL)$ via an edge map. The situation is summed up by the following commutative diagram, in which the vertical arrows are all inflation maps.\\
\[
\begin{tikzcd}
\HH^1(X,\LL_1)\otimes \HH^2(X,\LL_2) \arrow[r,"\cup"] & \HH^3(X,\LL_1\otimes\LL_2 ) \arrow[r,"\sim"]  & \HH^1(\GG,\HH^2(\bar X,\LL)) \\
~& ~& \HH^1(\SS^{(5)},H^{2,(5)}(M)) \arrow[u,"\sim"{anchor=north, rotate=90}]\\
~& \ker(\HH^3(G^{(4)},M)\to\HH^3(\bar G^{(4)},M))\arrow[uu]\arrow[r]& \HH^1(\SS^{(4)},\HH^2(\bar G^{(4)},M)) \arrow[u] \\
H^{1,(3)}(M_1)\otimes H^{2,(3)}(M_2)\arrow[r,"\cup"]\arrow[uuu,"\sim"{anchor=north, rotate=90}]& H^{3,(3)}(M)\arrow[u]& \\ ~
\HH^1(G^{(2)},M_1)\otimes \HH^2(G^{(2)},M_2)\arrow[r,"\cup"]\arrow[u]& \HH^3(G^{(2)},M)\arrow[u] & ~
\end{tikzcd}
\]
\end{proof}
\end{prop}

\vspace{-0.5cm}

\paragraph*{Acknowledgements} The author is grateful to the anonymous referee for the constructive comments and suggestions. He would also like to thank Quentin Remy for providing the aesthetically pleasing figure in \Cref{sec:expl}.

\bibliography{kp1_cup}
\bibliographystyle{plain}

\end{document}